\documentclass{clsmod}

\usepackage{latexsym,amssymb,graphicx}

\newtheorem{thm}{Theorem}[section]

\newtheorem{lem}[thm]{Lemma}

\newtheorem{cor}[thm]{Corollary}

\newnumbered{defn}[thm]{Definition}
\newnumbered{quest}[thm]{Question}
\newnumbered{conv}[thm]{Convention}
\newnumbered{rem}[thm]{Remark}
\newnumbered{exmp}[thm]{Example}
\newnumbered{ack}{Acknowledgment}

\newcommand{\hg}{{\widehat{\Gamma}}}
\newcommand{\hgo}{{\widehat{\Gamma}^{\rm op}}}

\newcommand{\lk}{\mbox{Lk}}
\newcommand{\st}{\mbox{St}}

\newcommand{\si}{{\mathcal S}}

\newcommand{\Z}{\ensuremath{{\mathbb{Z}}}}
\newcommand{\R}{\ensuremath{{\mathbb{R}}}}

\newcommand{\mH}{\ensuremath{{H}}}
\newcommand{\mA}{\ensuremath{{A}}}

\title[$\mH^*(\mA_\Gamma, \Z \mA_\Gamma)$]
{The cohomology of right angled Artin groups
with group ring coefficients}

\author{C. Jensen \and J. Meier}

\classno{20F36 (primary), 57M07 (secondary)}
     
\extraline{Jensen thanks the Louisiana Board of Regents for a Research
Competitiveness Subprogram grant.\\
\indent Meier thanks the American Mathematical Society for the
support of a Centennial Research Fellowship and Columbia University
for hosting him.\\
\indent Last revised on 11 December 2003.}

\begin{document}
\maketitle

\begin{abstract}
We give an explicit formula for the cohomology of a right
angled Artin group with group ring coefficients 
in terms of the cohomology of its
defining flag complex.
\end{abstract}

\section{Introduction}

Let $\Gamma$ be a finite simplicial graph and let
$\hg$ be the induced flag complex, i.e., the maximal
simplicial complex whose $1$-skeleton is $\Gamma$.
The associated \emph{right angled Artin group}
$A_\Gamma$ is the group presented by
\[
A_\Gamma = \langle V(\Gamma)~|~vw = wv\mbox{ if }
\{v,w\} \in E(\Gamma)\rangle\ .
\]
Because one can import topological properties
of the associated flag complex $\hg$ into the group 
$A_\Gamma$, these groups have provided important
examples of exotic behavior.  
(See for example  \cite{BeBr97}, \cite{BrMe01} and
\cite{LeNu03}.)
Here we refine the understanding of the end topology
of right angle Artin groups by giving an explicit formula for 
the cohomology of $A_\Gamma$ with
group ring coefficients in terms of the cohomology
of $\hg$ and links of simplices in $\hg$.

\begin{defn}  If $K$ is a simplicial complex let
$\si(K)$ denote the set of closed simplices --- including
the empty simplex --- in $K$.  The dimension of a simplex
is denoted $|\sigma|$; the link 
is denoted $\lk(\sigma)$;  the star of $\sigma$ is $\st(\sigma)$.
By definition $|\emptyset| = -1$ and $\lk(\emptyset) = K$.
\end{defn}

\begin{main*}
Let $\Gamma$ be a finite simplicial graph, let $\hg$ be
the associated flag complex and $A_\Gamma$ the associated
right angled Artin group.  As long as $\hg$ is not a 
single simplex, 
\[
H^*(A_\Gamma, \Z A_\Gamma) = \bigoplus_{\sigma \in \si(\hg)}
\left[ \bigoplus_{i=1}^\infty \overline{H}^{* - |\sigma|-2}
\left(\lk(\sigma)\right)
\right]\ .
\]
If $\hg$ is a single simplex then $A_\Gamma$ is free abelian
and $H^*(A_\Gamma, \Z A_\Gamma)$ is simply  $\Z$ in top
dimension.
\end{main*}

\begin{exmp}
Let $\hg$ be $\R P^2$.  
Then the reduced cohomology of $\lk(\emptyset) = \R P^2$ is concentrated
in dimension $2$ where it is $\Z_2$.  The link of any other simplex $\sigma$
is a $(1-|\sigma|)$-sphere hence its reduced cohomology is
concentrated in dimension $(1 - |\sigma|)$, where it is $\Z$.  Thus
$H^*(A_\Gamma, \Z A_\Gamma)$ is trivial except in dimension $3$ where
it is the sum of a countably generated free abelian group and a countable
sum of $\Z_2$'s.
\end{exmp}

There are at least two approaches to establishing the Main Theorem.
  One can modify the techniques of \cite{DaMe02}
that were developed 
for computing the cohomology of  Coxeter groups with
group ring coefficients --- as well as the cohomology with compact
supports of any locally finite building --- to compute this cohomology for
 right angled Artin groups.  In fact, the formula given in
the Main Theorem is quite similar to the formulas for cohomology
with compact supports of locally finite  buildings (Theorem~5.8 in
\cite{DaMe02}).  We  take a  more
efficient route, and 
 use the fact that right angled Artin groups
are commensurable with certain right angled Coxeter groups
\cite{DaJa02},
and appeal to the formula for the cohomology of a right angled
Coxeter group with group ring coefficients (\cite{Da98} or
\cite{DaMe02}). 

In the last section we explain how the formula of the
Main Theorem extends results of
\cite{BrMe01} on the end topology of right angled Artin groups.

\section{Background and Definitions}

One of the classical approaches to the study of 
asymptotic properties of a group $G$ is via its 
cohomology with $\Z G$-coefficients.  For example,  
from Proposition 7.5 and Exercise 4 of \cite{Br82},
if $G$ is a discrete group and $X$ is a contractible
$G$-complex with finite cell stabilizers and finite
quotient, then
$$H^*(G, \Z G) \cong H^*_c(X; \Z),$$
where $H^*_c(X; \Z)$ is the cohomology of $X$ with
compact supports.  In particular, one can take as
$X$ either of the classifying spaces $E G$
or $\underline{E} G$ 
provided they have finite quotients $BG$ or $\underline{B}G$ 
(cf. \cite{KrMi98}).
Cohomology  with group ring coefficients
determines 
the cohomological dimension of $G$ \cite[VIII.6.7]{Br82}:
If $G$ is of type FP then
\[
\mbox{cd } G = \hbox{ max}\{n:H^n(G, \Z G) \not = 0 \}.
\]
It is also closely related to connectivity at infinity
and duality properties as is described at the
end of the next section.

\begin{defn}
Right angled Artin groups admit CAT(0) $K(\pi,1)$s
formed as the union of tori.  
If $\Gamma$ is a finite
simplicial graph, 
let $K_\Gamma$ be the complex formed by joining tori
in the manner described by the flag complex $\hg$.  That is, for
each simplex $\sigma \subset \si(\hg)$, let $T_\sigma$
be the torus formed by identifying parallel faces of a unit
$(|\sigma|+1)$-cube.  (The torus $T_\emptyset$ is a single
vertex.)  The complex $K_\Gamma$ is then the union of
these tori, subject to $T_\sigma \cap T_{\sigma'} = T_{\sigma''}$
when $\sigma \cap \sigma' = \sigma''$ in $\hg$.  
For a proof that these $K_\Gamma$'s are CAT(0) classifying spaces,
see \cite{MV95}.  
We denote the  universal
cover of $K_\Gamma$ by $\widetilde{K}_\Gamma$.
\end{defn}

The complex $\widetilde{K}_\Gamma$ is also the Davis complex
for an appropriate  right angled Coxeter group.  Given
a finite simplicial graph $\Gamma$ the 
\emph{right angled Coxeter group} $C_\Gamma$ is the
quotient of $A_\Gamma$ formed by declaring that
each generator is an involution
\[
C_\Gamma = \langle V(\Gamma)~|~vw = wv\mbox{ if }
\{v,w\} \in E(\Gamma)\mbox{ and } v^2 = 1
\mbox{ for all }v\rangle\ .
\]

For a finite simplicial graph  $\Gamma$  let
$\Gamma'$ be the graph whose vertices are given
by $V(\Gamma) \times \{-1,1\}$ where
\begin{eqnarray*}
\{(v,\epsilon),(w, \epsilon)\} \in E(\Gamma')
   &\Leftrightarrow& \{v,w\} \in E(\Gamma) \\ 
\{(v,\epsilon), (w,- \epsilon)\} \in E(\Gamma')
   &\Leftrightarrow& v \ne w
\end{eqnarray*}
for $\epsilon = 1$ or $-1$. 


\begin{thm}[(Davis-Januszkiewicz \cite{DaJa02})]
\label{thm:commensurable}
The 
Artin group $A_\Gamma$ and the Coxeter group $C_{\Gamma'}$ are 
commensurable and in fact the complexes $\tilde K_\Gamma$ and the 
Davis complex 
for $C_{\Gamma'}$ are identical.
\end{thm}

\noindent (Because $\widetilde{K}_\Gamma$ is the Davis complex
for $C_{\Gamma'}$ we do not actually define the Davis complex for
a Coxeter group; see \cite{DaJa02} for a definition.)

One can now derive a formula for $H^*(A_\Gamma, \Z A_\Gamma)$
from known results in the literature.  Namely, because 
\begin{enumerate}
\item cohomology with group ring coefficients can be expressed
in terms of cohomology with compact supports of an 
$\underline{E} G$, and
\item $\widetilde{K}_\Gamma$ is both an $E A_\Gamma$ and
an $\underline{E} C_{\Gamma'}$, and
\item the cohomology of a Coxeter group with group ring
coefficients has been computed, and can be expressed
in terms of the cohomology of subcomplexes of links
of vertices in the Davis complex (\cite{Da98} or
\cite{DaMe02}), 
\end{enumerate}
we have the following formula for the cohomology  of $A_\Gamma$
with $\Z A_\Gamma$ coefficients. 

\begin{cor}\label{cor:firstformula}  Each $w \in C_{\Gamma'}$
has an associated simplex $\sigma(w) \in \si(\hg')$ such that
\[
H^*(A_\Gamma, \Z A_\Gamma) = H^*(C_{\Gamma'}, \Z C_{\Gamma'})
= H^*_c(\widetilde{K}_\Gamma; \Z) =  \bigoplus_{w \in C_{\Gamma'}}
\overline{H}^{*-1}\left(\hg' - \sigma(w)\right)\ .
\]
Each simplex $\sigma \in \si(\hg') \setminus \{\emptyset\}$
occurs countably many times in this sum, while 
$\sigma = \emptyset$ occurs exactly once.
\end{cor}

Although the formula above is correct, it obfuscates the connection
between $H^*(A_\Gamma, \Z A_\Gamma)$ and the cohomology
of the flag complex $\hg$.  
As a first step toward expressing the right hand side
in terms of the flag complex $\hg$, we give an alternate
description of the flag complex $\hg'$.

For each $v \in V(\Gamma)$ let $\hg_v$ be the full subcomplex
of $\hg$ 
induced by the vertices $V(\Gamma) \setminus \{v\}$. Thus
$\hg_v$ is a deformation retract of $\hg$ with the vertex
$v$ removed.

Let $(W,V(\Gamma))$ be the Coxeter system where $W$ is 
abelian and the generating set has been identified
with the vertices of the graph $\Gamma$.  Hence $W$ is simply 
\[
W = \underbrace{\Z_2 \times
\cdots \times \Z_2}_{|V(\Gamma)|\mbox{ copies}}\ .
\]
Let the $\hg_v$ be a set of mirrors related 
to this  Coxeter system
and form the associated $W$-complex in the
following manner.  For each $x \in \hg$
let $W_x$ be the subgroup of $W$ generated by the set of
$v \in V(\Gamma)$ such that $x$ belongs to $\hg_v$.
In other words, $W_x$ is generated by those $v$ such
that $x$ is not in
the open neighborhood of $v$ in $\hg$.  Define 
\[
L_\Gamma = W \times \hg/\sim
\]
where $(w,x) \sim (v,y)$ if and only if
$x = y\mbox{ and }w^{-1}v \in W_x$.

The complex  $\hg'$ shows up in the formula of
Corollary~\ref{cor:firstformula} because it is isomorphic to 
the link of any vertex  in $\widetilde{K}_\Gamma$.
One can find the following result in \cite{DaJa02}.

\begin{lem}
The complex $\hg'$ is isomorphic to $L_\Gamma$, and is 
isomorphic to 
the link of the vertex in $K_\Gamma$.
\end{lem}

If $\sigma \in \si(\hg)$ one can form a subcomplex
$L_{\sigma} \subset L_{\hg}$  by defining $W_\sigma$
to be the subgroup of $W$ 
generated by $\{v \in V(\Gamma)~|~v \not \in \sigma\}$, 
and forming $W_\sigma \times \hg/\sim$
where as before $(w,x) \sim (v,y)$ if and only if
$x = y\mbox{ and }w^{-1}v \in W_x$.  In particular,
if $\sigma = \emptyset$ (the empty simplex) then
$L_{\emptyset} = L_{\hg}$.

\begin{center}
\begin{figure}[hpbt]
\includegraphics[width=2.5in]{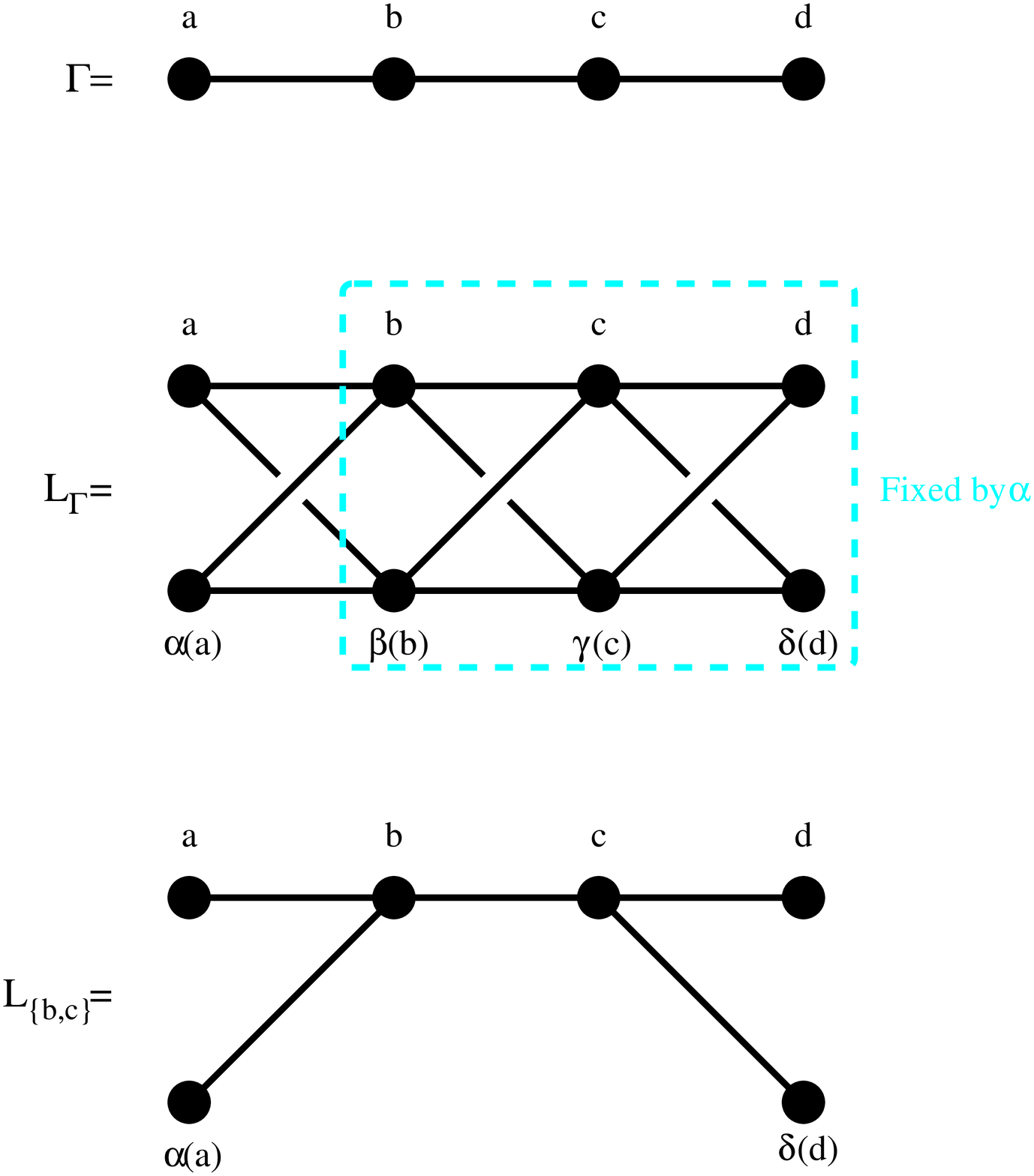}
\caption{A defining graph $\Gamma$, the associated
complex $L_\Gamma$, and subcomplex $L_{\{b,c\}}$ \label{fig:link}}
\end{figure}
\end{center}

\begin{exmp}\label{exmp:link}
Let $\Gamma = \hg$ be the simplicial arc indicated
in Figure~\ref{fig:link}.  The group $W$ is then generated
by four elements associated with the vertices.  Switching to
Greek letters we denote these generators as $\alpha, \beta,
\gamma$ and $\delta$, where  the mirror associated to
$\alpha$ is the subgraph induced by $\{b, c, d\}$, and
similarly for the other three generators.  
The complex $L_\Gamma$ is then as is indicated in
Figure~\ref{fig:link}.  The generator $\alpha$ acts on $L_\Gamma$
by exchanging the vertices labeled $a$ and $\alpha(a)$, 
and leaves all other vertices fixed.
Similarly $\beta$ exchanges $b$ and $\beta(b)$, fixing all other vertices, 
and so on.  Finally, if $\sigma = \{b,c\}$ then $L_\sigma$ is the bottom
complex in Figure~\ref{fig:link}.
\end{exmp}

For any $\sigma \in \si(\hg) \setminus \{\emptyset\}$ let
\[
\hg_{b(\sigma)} = \bigcup_{v \in \sigma^{(0)}} \hg_v\ ,
\]
so that $\hg_{b(\sigma)}$ is a deformation retract of
$\hg$ with the barycenter of $\sigma$ removed. 

\begin{lem}\label{lem:homologyformula}
The cohomology groups of $L_{\sigma}$ are given by
\[
\overline{H}^\ast(L_{\sigma}) = \bigoplus_{\tau \in \si(\hg - \sigma)} 
\overline{H}^\ast
\left(\hg, \hg_{b(\tau)}\right)\ ,
\]
where in a small abuse of notation we let $\si(\hg - \sigma)$ denote
all closed simplices of $\hg$ except those with non-empty intersection
with $\sigma$.
\end{lem} 

\begin{proof}
In \cite{Da87} Mike Davis gives a formula for the
homology of a complex on which a Coxeter group acts.
One can switch this to a formula for cohomology using
universal coefficients, or via  a minor
rewriting of Davis's original argument.  In our case
the formula is rather simple.  Since $W_\sigma$ is
abelian, each $w \in W_\sigma$ is determined by
the set of generators $S(w)$ that are necessary to
express $w$.  Temporarily following Davis's notation,
define
\[
\hg^{S(w)} = \bigcup_{v \in S(w)} \hg_v\ .
\]
(If $w = 1$ then $S(w) = \emptyset$ and so $\hg^{S(w)}$ is
empty as well.)  Davis's formula then gives
\[
\overline{H}^*(L_\sigma) \simeq \bigoplus_{w \in W_\sigma} 
      \overline{H}^*\left(\hg, \hg^{S(w)}\right)\ .
\]
This can be simplified.  If $S(w)$ is not the vertex set
of a simplex in $\hg$, then $\hg^{S(w)} = \hg$; if $S(w) = \sigma^{(0)}$
for some $\sigma \in \si(\hg)$, then $\hg^{S(w)} = \hg_{b(\sigma)}$.
Thus the formula above can be rewritten as
\[
\overline{H}^*(L_\sigma) \simeq \bigoplus_{\tau \in \si(\hg - \sigma)} 
       \overline{H}^*\left(\hg, \hg_{b(\tau)}\right)\ .
\]
\end{proof}

\section{Proof of the Main Theorem}

From Corollary~\ref{cor:firstformula} we know that
\[
H^*(A_\Gamma, \Z A_\Gamma) = \overline{H}^{*-1}\left(\hg'\right)
\bigoplus_{\sigma \in \si(\hg')\setminus \{\emptyset\}}
\left[ \bigoplus_{i=1}^\infty \overline{H}^{* - 1}
\left(\hg' - \sigma\right)
\right]
\]
where we know there are infinitely many copies of
$\overline{H}^{* - 1} \left(\hg' - \sigma\right)$ since
by its construction there are no non-trivial finite conjugacy classes
in $C_{\Gamma'}$. 
To arrive at our Main Theorem  we need a formula for
$\overline{H}^*(\hg' - \sigma)$ where $\sigma$ is any simplex
in $\si(\hg')$.
Thus our key lemma is:

\begin{lem}\label{lem:puncturedcoho}
Let $\sigma \in \si(\hg')$.  Then $\hg' - \sigma$
is homotopy equivalent to $L_\sigma$ and
\[
\overline{H}^*(L_\sigma) = 
\bigoplus_{\tau \in \si(\hg_\sigma)}
\overline{H}^{* - |\tau|-1}
\left(\lk(\tau)\right)\ .
\]
\end{lem} 

\begin{proof}
The complex $\hg$ embeds in $\hg'$ in a number of ways.  Let
the \emph{standard embedding} $\hg \hookrightarrow \hg'$ have image
the subcomplex induced by $\{(v,1)~|~v \in V(\Gamma)\}$.  Define
$\hgo$ to be the subcomplex induced by $\{(v,-1)~|~v \in V(\Gamma)\}$.
If $\sigma$ is a simplex in $\si(\hg')$ then $\sigma$ is defined
by a set of vertices in $\Gamma$ along with choices of $\pm 1$.
If 
\[\sigma \sim \{(a,1), (b,1), \ldots, 
       (c,1), (x,-1), \ldots , (y,-1), (z,-1)\}\]
 then the
automorphism $\alpha \beta \cdots \gamma$ takes $\sigma$ to the
simplex 
\[\sigma' \sim \{(a,-1), (b,-1), \ldots, 
       (c,-1), (x,-1), \ldots , (y,-1), (z,-1)\}.\]
  (Here we have used
the same convention on naming generators of $W$ as in
Example~\ref{exmp:link}.)  Thus  in
discussing the topology of $\hg' - \sigma$ for $\sigma \in \si(\hg')$,
we may without loss of generality assume 
$\sigma \subset \hgo \subset \hg'$.
But the space formed by 
removing the closed simplex $\sigma \subset \hgo$ from $\hg'$ 
deformation retracts onto the subcomplex formed
by making all possible reflections of $\hg$ that do not
involve  the generators of $W$ that correspond to
vertices of $\sigma$.  In other words, $\hg' - \sigma$
deformation retracts onto $L_\sigma$, which implies our first claim.

From Lemma~\ref{lem:homologyformula} we know
$\overline{H}^*(L_{\sigma}) = \bigoplus_{\tau \in \si(\hg - {\sigma})}
 \overline{H}^*\left(\hg, \hg_{b(\tau)}\right)$, thus it suffices to
establish
\[
\overline{H}^*\left(\hg, \hg_{b(\tau)}\right) \simeq 
\overline{H}^{* - |\tau|-1}
\left(\lk(\tau)\right)\ .
\]
First, if $\tau = \emptyset$, $\lk(\tau) = \hg$ and $|\tau| = -1$,
so we get $\overline{H}^{* - |\tau|-1}
\left(\lk(\tau)\right) = \overline{H}^*(\hg)$. If $\tau \ne \emptyset$ then
by excision, $\overline{H}^*(\hg, \hg_{b(\tau)}) \simeq
\overline{H}^*(\st(\tau), S^{|\tau|}\left[\lk(\tau)\right])$ where
$\st(\tau)$ is the closed star of $\tau$ and 
$S^i[\cdot]$ denotes the $i^{th}$ suspension.
Because the star $\st(\tau)$ is contractible,
the long exact sequence in cohomology shows
\[
\overline{H}^*(\st(\tau), S^{|\tau|}\left[\lk(\tau)\right]) 
= \overline{H}^{*-1}(S^{|\tau|}\left[\lk(\tau)\right])\ .
\]
But the cohomology of a suspension is just a shifted
copy of the cohomology of the original complex
\[
\overline{H}^{*-1}(S^{|\tau|}\left[\lk(\tau)\right])
= \overline{H}^{* - |\tau|-1}\left(\lk(\tau)\right)
\]
and the result follows.
\end{proof}

\begin{exmp}
In Example~\ref{exmp:link} we considered $\Gamma = \hg = $ a simplicial
arc, and two associated complexes, $L_\Gamma$ and 
$L_{\{b,c\}}$.   (The first claim of Lemma~\ref{lem:puncturedcoho}
states that $L_{\{b,c\}}$ is homotopy equivalent to $L_\Gamma$ with
the closed edge $\{\beta(b), \gamma(c)\}$ removed.)
 The formula of Lemma~\ref{lem:puncturedcoho} says,
for example, that
\[
\overline{H}^1(L_\Gamma) = \bigoplus_{\sigma \in \si(\hg)} 
    \overline{H}^{1 - |\sigma| - 1}\left(\lk(\sigma)\right)
    = \bigoplus_{\sigma \in \si(\hg)}\overline{H}^{-|\sigma|}(\lk(\sigma))\ .
\]
This then becomes
\begin{eqnarray*}
H^1(L_\Gamma) &=& \overline{H}^1(\lk(\emptyset)) 
 \oplus \overline{H}^0(\lk(a)) 
 \oplus \overline{H}^0(\lk(b))  \oplus \overline{H}^0(\lk(c))  
 \oplus \overline{H}^0(\lk(d))   \\
  &\ & \oplus \overline{H}^{-1}(\lk(\{a,b\}) 
 \oplus \overline{H}^{-1}(\lk(\{b,c\})
 \oplus \overline{H}^{-1}(\lk(\{c,d\}) \\
  &=& \overline{H}^1(\hg) \oplus \left(\overline{H}^0(\bullet)\right)^2 
   \oplus \left(\overline{H}^0(\bullet\ \  \bullet)\right)^2
    \oplus \left(\overline{H}^{-1}(\emptyset)\right)^3 = \Z^5\ ,
\end{eqnarray*}
using the convention that $\overline{H}^{-1}(\emptyset) = \Z$.

In the case of $L_{\{b,c\}}$ one drops all the terms involving
$b$ or $c$, which are precisely the non-trivial terms above, 
hence $H^1(L_{\{b,c\}}) = 0$.
\end{exmp}

We can now prove our Main Theorem.

\begin{thm}
Let $\Gamma$ be a finite simplicial graph, let $\hg$ be
the associated flag complex and $A_\Gamma$ the associated
right angled Artin group.  As long as $\hg$ is not a 
single simplex, 
\[
H^*(A_\Gamma, \Z A_\Gamma) = \bigoplus_{\sigma \in \si(\hg)}
\left[ \bigoplus_{i=1}^\infty \overline{H}^{* - |\sigma|-2}
\left(\lk(\sigma)\right)
\right]\ .
\]
If $\hg$ is a single simplex then $A_\Gamma$ is free abelian
and $H^*(A_\Gamma, \Z A_\Gamma)$ is simply $\Z$ in top
dimension.
\end{thm}

\begin{proof}
From Corollary~\ref{cor:firstformula} we have
\[
H^*(A_\Gamma, \Z A_\Gamma) = \bigoplus_{w \in C_{\Gamma'}}
\overline{H}^{*-1}\left(\hg' - \sigma(w)\right)\ .
\]
By Lemma~\ref{lem:puncturedcoho} this gives
\[
H^*(A_\Gamma, \Z A_\Gamma) = \bigoplus_{w \in C_{\Gamma'}}
\left[ \bigoplus_{\tau \in \si(\hg - \sigma(w))}
\overline{H}^{*-|\tau| - 2}\left(\lk(\tau)\right)
\right]\ .
\]
If $\hg$ is not a single simplex, then each $\tau \in \si(\hg)$
will show up in the product inside the square brackets 
for infinitely many $w \in C_{\Gamma'}$, and the formula 
in the theorem follows.  

On the other hand, if $\hg$ is a single
simplex $\sigma$, then $\sigma$ only occurs in the summand
corresponding to $1 \in C_{\Gamma'}$.  All other simplices
occur infinitely often, but if $\tau \ne \sigma$ then
 $\lk(\tau)$ is contractible, and $\overline{H}^*(\lk(\tau))$ is
zero.  Thus $H^*(A_\Gamma, \Z A_\Gamma) 
= H^{* - |\sigma| - 2}(\emptyset)$, consistent with
the fact that $A_\Gamma = \Z^{|\sigma| + 1}$.
\end{proof}

As was alluded to in the previous section,
cohomology with group ring coefficients is
closely related to asymptotic properties.  
A group $G$ that admits a finite $K(G,1)$
is \emph{$n$-acyclic at infinity} if roughly
speaking, complements of compact sets in the
universal cover have trivial homology through
dimension $n$ (see \cite{DaMe02} for a precise definition.)
It was from this perspective that Brady and Meier determined
when a right angled Artin group was $n$-acyclic
at infinity.  Their approach was via a combinatorial
 Morse theory argument using the $K_\Gamma$ complexes.
However, there is an algebraic characterization that says 
a group $G$ is $n$-acyclic at infinity if and only
if $H^i(G,\Z G) = 0$ for $i \le n+1$ and 
$H^{n+2}(G, \Z G)$ is torsion-free (see \cite{GeMi86}).
The group $G$ is an \emph{$n$-dimensional duality group}
if there is a dualizing module $D$ such that
$H_i(G,M) \simeq H^{n-i}(G,M \otimes D)$ for all $i$ and
all $G$-modules $M$.  This too can be recast in terms
of cohomology with group ring coefficients:
$G$ is an $n$-dimensional duality group
if its cohomology with group ring coefficients is
torsion-free and concentrated
in dimension $n$ \cite{Bi81}.
Thus our Main Theorem  implies three results
of \cite{BrMe01}.  It is important to remember 
that $\emptyset \in \si(\hg)$, and
the formal dimension of $\emptyset$ is $-1$.

\begin{cor}[(Prop.~4.1 in \cite{BrMe01})]\label{cor:acyclicity}
A right angled Artin group $A_\Gamma$ is $n$-acyclic at
infinity if and only if for all $\sigma \in \si(\hg)$,
$\lk(\sigma)$  is $(n - |\sigma| - 1)$-acyclic.
\end{cor}

\begin{proof}
Since $\lk(\sigma)$ is $(n-|\sigma|-1)$-acyclic it
follows by universal coefficients
 that its cohomology is trivial up to dimension $n-|\sigma|-1$
and that $\overline{H}^{n-|\sigma|}(\lk(\sigma))$ is torsion-free.
Thus the formula of the Main Theorem implies that 
$H^i(A_\Gamma, \Z A_\Gamma)$ is zero for $i \le n+1$ and
$H^{n+2}(A_\Gamma, \Z A_\Gamma)$ is torsion-free.
\end{proof}

\begin{cor}[(Theorem C in \cite{BrMe01})]\label{cor:duality}
A right angled Artin group $A_\Gamma$ is a duality group
if and only if $\hg$ is Cohen-Macaulay.
\end{cor}

\begin{proof}
 A simplicial complex $K$ is \emph{Cohen-Macaulay} if
for any simplex $\sigma \in \si(K)$, the cohomology of
$\lk(\sigma)$ is concentrated in top dimension (and is 
torsion free).  It follows from the formula of the Main Theorem
that $H^*(A_\Gamma, \Z A_\Gamma)$ is torsion free and concentrated
in top dimension if and only if $\hg$ is Cohen-Macaulay.
\end{proof}

Recall that an $n$-dimensional
duality group is called  a \emph{Poincar\'e} duality group if and only
if $H^n(G, \Z G) = \Z$ \cite{BiEck73}.
After the statement of Theorem~C in \cite{BrMe01} it was
remarked that a Theorem of Strebel combined with Theorem~C implies that
 a right angled Artin group $A_\Gamma$ is a Poincar\'e duality
group if and only if $A_\Gamma$ is free abelian.  
This characterization follows directly
from the formula in our Main Theorem.

\begin{cor}\label{cor:poincare}
A right angled Artin group $A_\Gamma$ is a Poincar\'e
duality group if and only if it is free abelian.
\end{cor}

\begin{proof}
The Main Theorem implies that $H^n (A_\Gamma, \Z A_\Gamma)$ is
not finitely generated --- in particular it is not equal to $\Z$ 
 ---  unless $\hg$ is a simplex and hence $A_\Gamma$
is free abelian.
\end{proof}


\affiliationone{ Craig A.~Jensen\\
               Department of Mathematics\\
               University of New Orleans\\
               New Orleans, LA 70148
      \email{jensen@math.uno.edu}}
\affiliationtwo{ John Meier\\
               Department of Mathematics\\
               Lafayette College\\
               Easton, PA 18042
      \email{meierj@lafayette.edu}}

\end{document}